\documentclass[11pt, leqno]{article}

\usepackage{latexsym,graphicx,times,amsmath,amsfonts}

\newtheorem{theorem}{Theorem}
\newtheorem{lemma}[theorem]{Lemma}
\newtheorem{corollary}[theorem]{Corollary}

\newenvironment{Theorem}
	{\begin{theorem}\sl}
	{\end{theorem}}
\newenvironment{Lemma}
	{\begin{lemma}\sl}
	{\end{lemma}}
\newenvironment{Corollary}
	{\begin{corollary}\sl}
	{\end{corollary}}

\def\bs{\boldsymbol}

\def\TV{\mathrm{TV}}

\def\N{\mathbb{N}}
\def\Z{\mathbb{Z}}
\def\R{\mathbb{R}}
\def\Pr{\mathbb{P}}

\begin{document}

\addtolength{\baselineskip}{+.49\baselineskip}

\title{\bf Explicit Bounds for the Approximation Error in Benford's Law}
\author{Lutz D\"{u}mbgen and Christoph Leuenberger\\
	\it{University of Berne and Ecole d'Ing\'{e}nieurs de Fribourg}}
\date{June 2007, revised December 2007}
\maketitle

\vfill

\paragraph{Abstract.}
Benford's law states that for many random variables $X > 0$ its leading digit $D = D(X)$ satisfies approximately the equation $\Pr(D = d) = \log_{10}(1 + 1/d)$ for $d = 1,2,\ldots,9$. This phenomenon follows from another, maybe more intuitive fact, applied to $Y := \log_{10}X$: For many real random variables $Y$, the remainder $U := Y - \lfloor Y\rfloor$ is approximately uniformly distributed on $[0,1)$. The present paper provides new explicit bounds for the latter approximation in terms of the total variation of the density of $Y$ or some derivative of it. These bounds are an interesting and powerful alternative to Fourier methods. As a by-product we obtain explicit bounds for the approximation error in Benford's law.
\vfill


\newpage

\section{Introduction}
\label{Introduction}

The First Digit Law is the empirical observation that in many tables of numerical data the leading significant digits are not uniformly distributed as one might suspect at first. The following law was first postulated by Simon Newcomb (1881):
$$
	\mbox{Prob}(\mbox{leading digit} = d) = \log_{10}(1 + 1/d)
$$
for $d=1,\dots,9$. Since the rediscovery of this distribution by physicist Frank Benford (1938), an abundance of additional empirical evidence and various extensions have appeared, see Raimi (1976) and Hill (1995) for a review. Examples for ``Benford's law'' are one-day returns on stock market indices, the population sizes of U.S.\ counties, or stream flow data (Miller and Nigrini 2007). An interesting application of this law is the detection of accounting fraud (see Nigrini, 1996). Numerous number sequences (e.g.\ Fibonacci's sequence) are known to follow Benford's law exactly, see Diaconis (1977), Knuth (1969) and Jolissaint (2005).

An elegant way to explain and extend Benford's law is to consider a random variable $X > 0$ and its expansion with integer base $b \ge 2$. That means, $X = M \cdot b^Z$ for some integer $Z$ and some number $M \in [1,B)$, called the mantissa of $X$. The latter may be written as $M = \sum_{i=0}^\infty D_i \cdot b^{-i}$ with digits $D_i \in \{0,1,\ldots,b-1\}$. This expansion is unique if we require that $D_i \ne b-1$ for infinitely many indices $i$, and this entails that $D_0 \ge 1$. Then the $\ell + 1$ leading digits of $X$ are equal to $d_0, \ldots, d_\ell \in \{0,1,\ldots,b-1\}$ with $d_0 \ge 1$ if, and only if,
\begin{equation}
	\bs{d} \ \le \ M \ < \ \bs{d} + b^{-\ell}
	\quad\text{with}\quad
	\bs{d} \ := \ \sum_{i=0}^\ell d_i \cdot b^{-i} .
	\label{Benford-event}
\end{equation}
In terms of $Y := \log_b(X)$ and
$$
	U \ := \ Y - \lfloor Y\rfloor \ = \ \log_b(M)
$$
one may express the probability of (\ref{Benford-event}) as
\begin{equation}
	\Pr \bigl( \log_b(\bs{d}) \le U < \log_b(\bs{d} + b^{-\ell}) \bigr) .
	\label{Benford-probability}
\end{equation}
If the distribution of $Y$ is sufficiently ``diffuse'', one would expect the distribution of $U$ being approximately uniform on $[0,1)$, so that (\ref{Benford-probability}) is approximately equal to
$$
	\log_b(\bs{d} + b^{-\ell}) - \log_b(\bs{d})
	\ = \ \log_b(1 + b^{-\ell}/\bs{d}) .
$$

Hill (1995) stated the problem of finding distributions satisfying Benford's law exactly. Of course, a sufficient condition would be $U$ being uniformly distributed on $[0,1)$. Leemis et al.\ (2000) tested the conformance of several survival distributions to Benford's law using computer simulations. The special case of exponentially distributed random variables was studied by Engel and Leuenberger (2003): Such random variables satisfy the first digit law only approximatively, but precise estimates can be given; see also Miller and Nigrini (2006) for an alternative proof and extensions. Hill and Schuerger (2005) study the regularity of digits of random variables in detail.

In general, uniformity of $U$ isn't satisfied exactly but only approximately. Here is one typical result: Let $Y = \sigma Y_o$ for some random variable $Y_o$ with Lebesgue density $f_o$ on the real line. Then
$$
	\sup_{B \in \mathrm{Borel}([0,1))} \, \bigl| \Pr(U \in B) - \mathrm{Leb}(B) \bigr| \ \to \ 0
	\quad\text{as} \ \sigma \to \infty .
$$
This particular and similar results are typically derived via Fourier methods; see, for instance, Pinkham (1961) or Kontorovich and Miller (2005).

The purpose of the present paper is to study approximate uniformity of the remainder $U$ in more detail. In particular we refine and extend an inequality of Pinkham (1961). Section~\ref{Distribution of U} provides the density and distribution function of $U$ in case of the random variable $Y$ having Lebesgue density $f$. In case of $f$ having finite total variation or, alternatively, $f$ being $k \ge 1$ times differentiable with $k$-th derivative having finite total variation, the deviation of $\mathcal{L}(U)$ (i.e.\ the distribution of $U$) from $\mathrm{Unif}[0,1)$ may be bounded explicitly in several ways. Since any density may be approximated in $L^1(\R)$ by densities with finite total variation, our approach is no less general than the Fourier method. Section~\ref{Applications} contains some specific applications of our bounds. For instance, we show that in case of $Y$ being normally distributed with variance one or more, the distribution of the remainder $U$ is \textsl{very} close to the uniform distribution on $[0,1)$.

\section{On the distribution of the remainder $U$}
\label{Distribution of U}

Throughout this section we assume that $Y$ is a real random variable with c.d.f.\ $F$ and Lebesgue density $f$.

\subsection{The c.d.f.\ and density of $U$}

For any Borel set $B \subset [0,1)$,
$$
	\Pr(U \in B) \ = \ \sum_{z \in \Z} \Pr(Y \in z + B) .
$$
This entails that the c.d.f.\ $G$ of $U$ is given by
$$
	G(x) := \Pr(U \le x) \ = \ \sum_{z \in \Z} (F(z+x) - F(z))
	\quad\text{for} \ 0 \le x \le 1 .
$$
The corresponding density $g$ is given by
$$
	g(x) \ := \ \sum_{z \in \Z} f(z + x) .
$$
Note that the latter equation defines a periodic function $g : \R \to [0,\infty]$, i.e.\ $g(x+z) = g(x)$ for arbitrary $x \in \R$ and $z \in \Z$. Strictly speaking, a density of $U$ is given by $1\{0 \le x < 1\} g(x)$.

\subsection{Total variation of functions}

Let us recall the definition of total variation (cf.\ Royden 1988, Chapter~5): For any interval $\mathbb{J} \subset \R$ and a function $h : \mathbb{J} \to \R$, the total variation of $h$ on $\mathbb{J}$ is defined as
$$
	\TV(h, \mathbb{J}) \ := \ \sup \Bigl\{ \sum_{i=1}^m \bigl| h(t_i) - h(t_{i-1}) \bigr| \, : \
		m \in \N; \, t_0 < t_1 < \cdots < t_m; \, t_0,\ldots,t_m \in \mathbb{J} \Bigr\} .
$$
In case of $\mathbb{J} = \R$ we just write $\TV(h) := \TV(h, \R)$. If $h$ is absolutely continuous with derivative $h'$ in $L_{\rm loc}^1(\R)$, then
$$
	\TV(h) \ = \ \int_{\R} |h'(x)| \, dx .
$$
An important special case are unimodal probability densities $f$ on the real line, i.e.\ $f$ is non-decreasing on $(-\infty,\mu]$ and non-increasing on $[\mu,\infty)$ for some real number $\mu$. Here $\TV(f) = 2 f(\mu)$.

\subsection{Main results}

We shall quantify the distance between $\mathcal{L}(U)$ and $\mathrm{Unif}[0,1)$ by means of the range of $g$,
$$
	\mathrm{R}(g)
	\ := \ \sup_{x, y \in \R} \bigl| g(y) - g(x) \bigr|
	\ \ge \ \sup_{u \in [0,1]} |g(u) - 1| .
$$
The latter inequality follows from $\sup_{x\in\R} g(x) \ge \int_0^1 g(x) \, dx = 1 \ge \inf_{x\in\R} g(x)$. In addition we shall consider the Kuiper distance between $\mathcal{L}(U)$ and $\mathrm{Unif}[0,1)$,
$$
	\mathrm{KD}(G)
	\ := \ \sup_{0 \le x < y \le 1} \bigl| G(y) - G(x) - (y - x) \bigr|
	\ = \ \sup_{0 \le x < y \le 1} \bigl| \Pr(x \le U < y) - (y - x) \bigr| ,
$$
and the maximal relative approximation error,
$$
	\mathrm{MRAE}(G)
	\ := \ \sup_{0 \le x < y \le 1} \Bigl| \frac{G(y) - G(x)}{y - x} - 1 \Bigr| .
$$
Expression (\ref{Benford-probability}) shows that these distance measures are canonical in connection with Benfords law. Note that $\mathrm{KD}(G)$ is bounded from below by the more standard Kolmogorov-Smirnov distance,
$$
	\sup_{x \in [0,1]} |G(x) - x| ,
$$
and it is not greater than twice the Kolmogorov-Smirnov distance.

\begin{Theorem}
\label{thm: TV(f0)}
Suppose that $\TV(f) < \infty$. Then $g$ is real-valued with
$$
	\TV(g, [0,1]) \ \le \ \TV(f)
	\quad\text{and}\quad
	\mathrm{R}(g) \ \le \ \TV(f)/2 .
$$
\end{Theorem}

\paragraph{Remark.}
The inequalities in Theorem~\ref{thm: TV(f0)} are sharp in the sense that for each number $\tau > 0$ there exists a density $f$ such that the corresponding density $g$ satisfies
\begin{equation}
	\TV(g, [0,1]) \ = \ \TV(f) \ = \ 2\tau
	\quad\text{and}\quad
	\max_{0 \le x < y \le 1} \, \bigl| g(x) - g(y) \bigr| \ = \ \tau .
	\label{TV-equalities}
\end{equation}
A simple example, mentioned by the referee, is the uniform density $f(x) = 1\{0 < x < \tau\}/\tau$. Writing $\tau = m + a$ for some integer $m \ge 0$ and $a \in (0,1]$, one can easily verify that
$$
	g(x) \ = \ m/\tau + 1\{0 < x < a\}/\tau ,
$$
and this entails (\ref{TV-equalities}).

\begin{figure}[h]
\includegraphics[width=7.4cm]{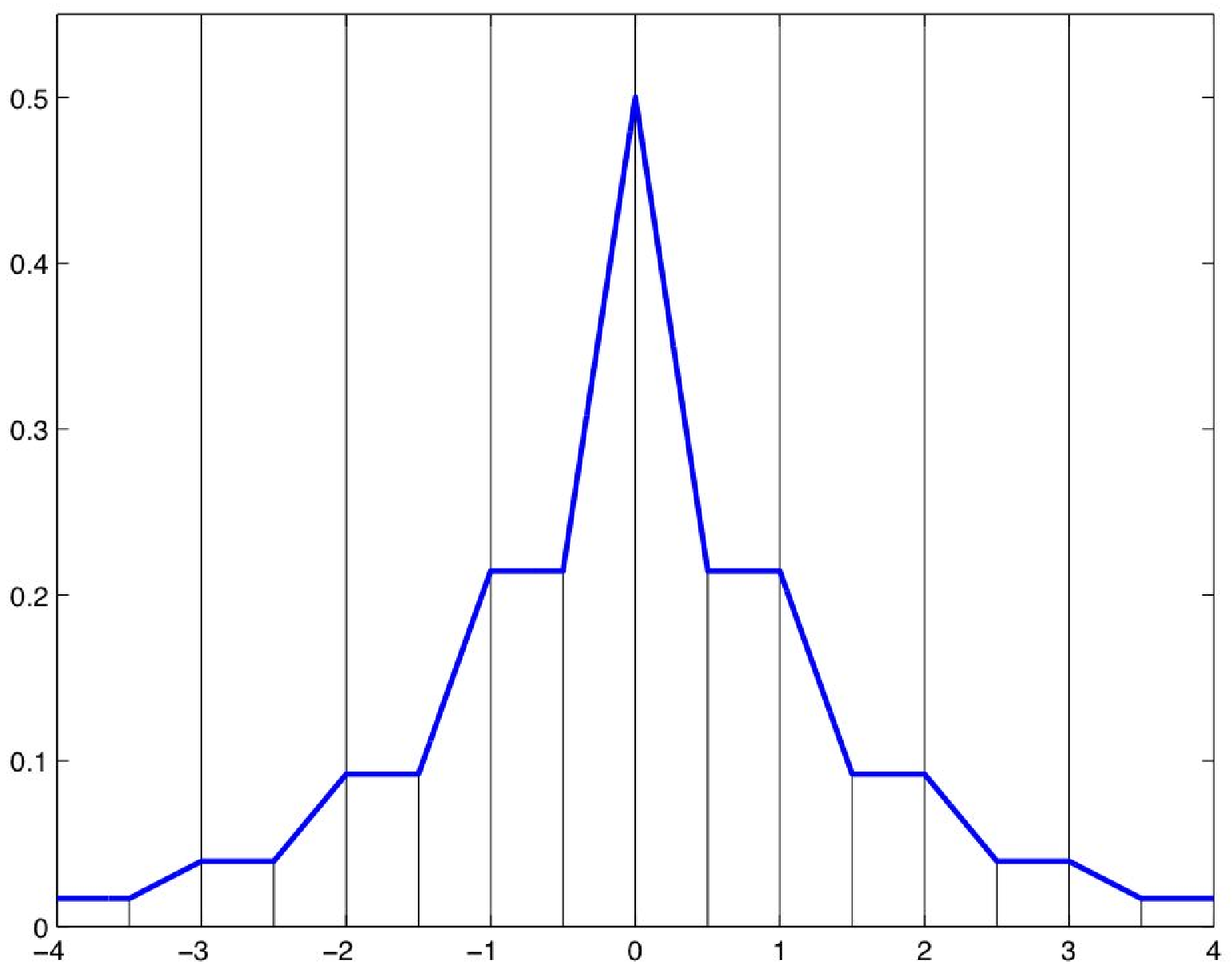}
\hfill
\includegraphics[width=7.4cm]{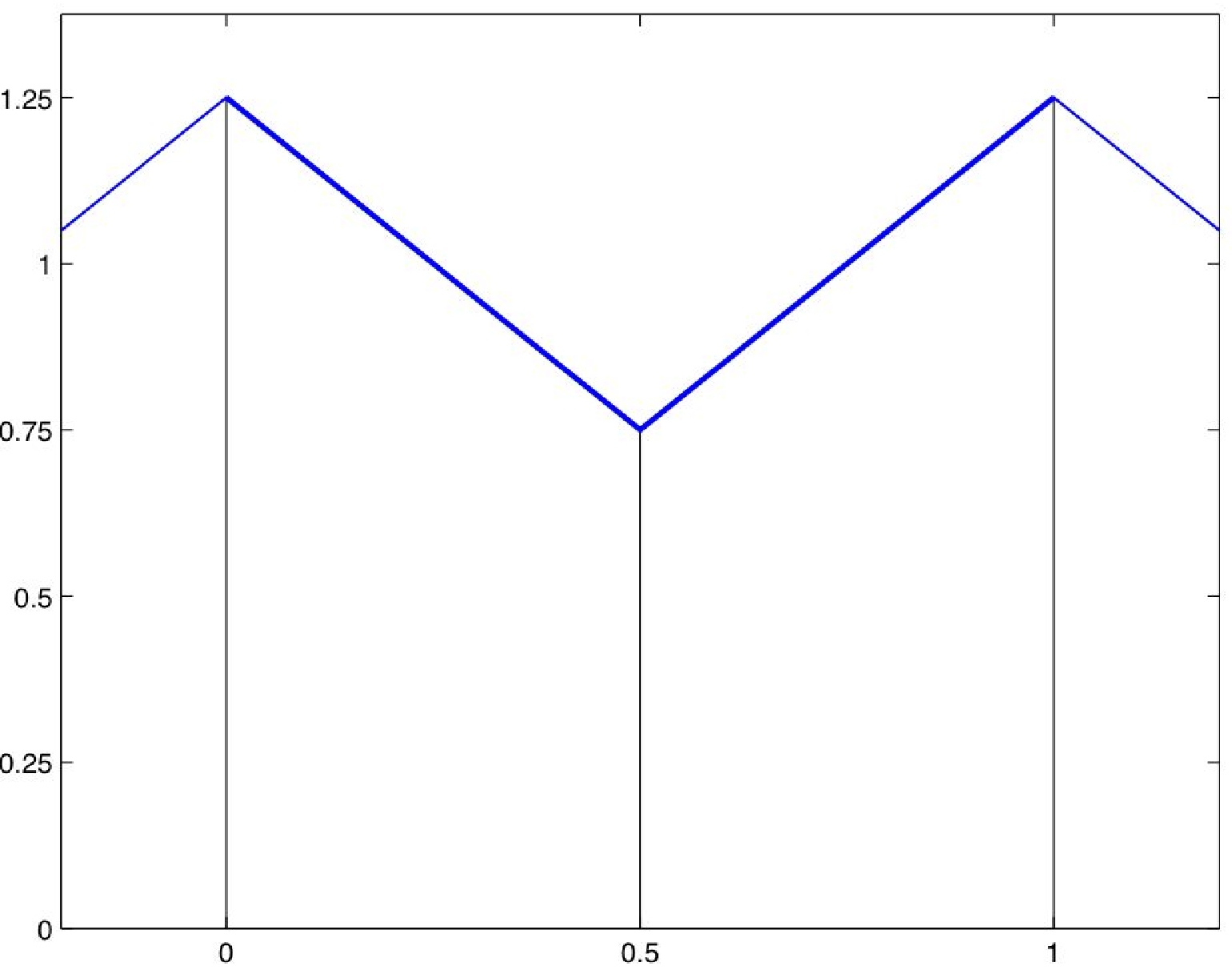}
\caption{A density $f$ (left) and the corresponding $g$ (right) such that $\TV(f) = \TV(g)$.}
\label{fig: WorstCase}
\end{figure}

Here is another example with continuous densities $f$ and $g$: For given $\tau > 0$ consider a continuous, even density $f$ with $f(0) = \tau$ such that for all integers $z \ge 0$,
$$
	f \ \text{is} \ \begin{cases}
		\text{linear and non-increasing on} \ [z, z+1/2] , \\
		\text{constant on} \ [z+1/2, z+1] .
	\end{cases}
$$
Then $f$ is unimodal with mode at zero, whence $\TV(f) = 2 f(0) = 2 \tau$. Moreover, one verifies easily that $g$ is linear and decreasing on $[0,1/2]$ and linear and increasing on $[1/2,1]$ with $g(0) - g(1/2) = \tau$. Thus $\TV(g,[0,1]) = 2\tau$ as well. Figure~\ref{fig: WorstCase} illustrates this construction. The left panel shows (parts of) an even density $f$ with $f(0) = 0.5 = \TV(f)/2$, and the resulting function $g$ with $\TV(g,[0,1]) = \TV(f) = g(1) - g(0.5)$.

As a corollary to Theorem~\ref{thm: TV(f0)} we obtain a refinement of the inequality
$$
	\sup_{0 \le x \le 1} |G(x) - x| \le \TV(f)/6
$$
which was obtained by Pinkham (1961, corollary to Theorem~2) via Fourier techniques:

\begin{Corollary}
\label{cor: TV(f0)}
Under the conditions of Theorem~\ref{thm: TV(f0)}, for $0 \le x < y \le 1$,
$$
	\bigl| G(y) - G(x) - (y - x) \bigr|
	\ \le \ (y - x)(1 - (y-x)) \TV(f) / 2 .
$$
In particular,
$$
	\mathrm{KD}(G) \ \le \ \TV(f)/8
	\quad\text{and}\quad
	\mathrm{MRAE}(G) \ \le \ \TV(f)/2 .
$$
\end{Corollary}

The previous results are for the case of $\TV(f)$ being finite. Next we consider smooth densities $f$. A function $h$ on the real line is called $k \ge 1$ times absolutely continuous if $h \in \mathcal{C}^{k-1}(\R)$, and if its derivative $h^{(k-1)}$ is absolutely continuous. With $h^{(k)}$ we denote some version of the derivative of $h^{(k-1)}$ in $L_{\rm loc}^1(\R)$.

\begin{Theorem}
\label{thm: TV(fk)}
Suppose that $f$ is $k \ge 1$ times absolutely continuous such that $\TV(f^{(k)}) < \infty$ for some version of $f^{(k)}$. Then $g$ is Lipschitz-continuous on $\R$. Precisely, for $x,y \in \R$ with $|x - y| \le 1$,
$$
	\bigl| g(x) - g(y) \bigr| \ \le \ |x - y| (1 - |x - y|) \frac{\TV(f^{(k)})}{2 \cdot 6^{k-1}}
	\ \le \ \frac{\TV(f^{(k)})}{8 \cdot 6^{k-1}} .
$$
\end{Theorem}

\begin{Corollary}
\label{cor: TV(fk)}
Under the conditions of Theorem~\ref{thm: TV(fk)}, for $0 \le x < y \le 1$,
$$
	\bigl| G(y) - G(x) - (y - x) \bigr|
	\ \le \ (y - x)(1 - (y-x)) \frac{\TV(f^{(k)})}{2 \cdot 6^k} .
$$
In particular,
$$
	\mathrm{KD}(G) \ \le \ \frac{\TV(f^{(k)})}{8 \cdot 6^k}
	\quad\text{and}\quad
	\mathrm{MRAE}(G) \ \le \ \frac{\TV(f^{(k)})}{2 \cdot 6^k} .
$$
\end{Corollary}

Finally, let us note that Theorem~\ref{thm: TV(f0)} entails a short proof of the qualitative result mentioned in the introduction:

\begin{Corollary}
\label{cor: f0}
Let $Y = \mu + \sigma Y_o$ for some $\mu \in \R$, $\sigma > 0$ and a random variable $Y_o$ with density $f_o$, i.e. $f(x) = f_o((x - \mu)/\sigma)/\sigma$. Then
$$
	\int_0^1 |g(x) - 1| \, dx \ \to \ 0
	\quad\text{as} \ \sigma \to \infty, \ \text{uniformly in} \ \mu .
$$
\end{Corollary}

\section{Some applications}
\label{Applications}

We start with a general remark on location-scale families. Let $f_o$ be a probability density on the real line such that $\TV(f_o^{(k)}) < \infty$ for some integer $k \ge 0$. For $\mu \in \R$ and $\sigma > 0$ let
$$
	f(x) = f_{\mu,\sigma}(x) \ := \ \sigma^{-1} f \bigl(\sigma^{-1} (x - \mu) \bigr) .
$$
Then one verifies easily that
$$
	\TV(f^{(k)}) \ = \ \TV(f_o^{(k)}) / \sigma^{k+1} .
$$

\subsection{Normal and log-normal distributions}
\label{subsec: Benford-Gauss}

For $\phi(x) := (2\pi)^{-1/2}\exp(- x^2/2)$, elementary calculations reveal that
\begin{eqnarray*}
	\TV(\phi)
	& = & 2 \phi(0) \ \approx \ 0.7979 , \\
	\TV(\phi^{(1)})
	& = & 4 \phi(1) \ \approx \ 0.9679 , \\
	\TV(\phi^{(2)})
	& = & 8 \phi(\sqrt{3}) + 2 \phi(0) \ \approx \ 1.5100 .
\end{eqnarray*}
In general,
$$
	\phi^{(k)}(x) \ = \ H_k(x) \phi(x)
$$
with the Hermite type polynomial
$$
	H_k(x) \ = \ \exp(x^2/2) \, \frac{d^k}{dx^k} \, \exp(- x^2/2)
$$
of degree $k$. Via partial integration and induction one may show that
$$
	\int H_j(x) H_k(x) \phi(x) \, dx \ = \ 1\{j = k\} k!
$$
for arbitrary integers $j, k \ge 0$ (cf.\ Abramowitz and Stegun 1964). Hence the Cauchy-Schwarz inequality entails that
\begin{eqnarray*}
	\TV(\phi^{(k)})
	& = & \int |\phi^{(k+1)}(x)| \, dx \\
	& = & \int |H_{k+1}(x)| \phi(x) \, dx \\
	& \le & \Bigl( \int H_{k+1}(x)^2 \phi(x) \, dx \Bigr)^{1/2} \\
	& = & \sqrt{ (k+1)! } .
\end{eqnarray*}
These bounds yield the following results:

\begin{Theorem}
\label{thm: Benford-Gauss}
Let $f(x) = f_{\mu,\sigma}(x) = \phi((x - \mu)/\sigma)/\sigma$ for $\mu \in \R$ and $\sigma \ge 1/6$. Then the corresponding functions $g = g_{\mu,\sigma}$ and $G = G_{\mu,\sigma}$ satisfy the inequalities
$$
	\begin{array}{rcr}
	\mathrm{R}(g_{\mu,\sigma})
	& \le & 4.5 \cdot  h \bigl( \lfloor 36 \sigma^2 \rfloor \bigr) , \\[1ex]
	\mathrm{KD}(G_{\mu,\sigma})
	& \le & 0.75 \cdot  h \bigl( \lfloor 36 \sigma^2 \rfloor \bigr) , \\[1ex]
	\mathrm{MRAE}(G_{\mu,\sigma})
	& \le & 3 \cdot  h \bigl( \lfloor 36 \sigma^2 \rfloor \bigr) ,
	\end{array}
$$
where $h(m) := \sqrt{m! / m^m}$ for integers $m \ge 1$.
\end{Theorem}

It follows from Stirling's formula that $h(m) = c_m m^{1/4} e^{-m/2}$ with $\lim_{m \to \infty} c_m = (2\pi)^{1/4}$. In particular,
$$
	\lim_{m \to \infty} \frac{\log h(m)}{m} \ = \ - \, \frac{1}{2} ,
$$
so the bounds in Theorem~\ref{thm: Benford-Gauss} decrease exponentially in $\sigma^2$. For $\sigma = 1$ we obtain already the remarkable bounds
$$
	\begin{array}{rcrcl}
	\mathrm{R}(g)
	& \le & 4.5 \cdot h(36) & \approx & 2.661 \cdot 10^{-7} , \\[1ex]
	\mathrm{KD}(G)
	& \le & 0.75 \cdot h(36) & \approx & 4.435 \cdot 10^{-8} , \\[1ex]
	\mathrm{MRAE}(G)
	& \le & 3 \cdot h(36) & \approx & 1.774 \cdot 10^{-7}
	\end{array}
$$
for all normal densities $f$ with standard deviation at least one.

\begin{Corollary}
\label{cor: Benford-Gauss}
For an integer base $b \ge 2$ let $X = b^Y$ for some random variable $Y \sim \mathcal{N}(\mu, \sigma^2)$ with $\sigma \ge 1/6$. Then the leading digits $D_0, D_1, D_2, \ldots$ of $X$ satisfy the following inequalities: For arbitrary digits $d_0, d_1, d_2, \ldots \in \{0,1,\ldots,b-1\}$ with $d_0 \ge 1$ and integers $\ell \ge 0$,
$$
	\left| \frac{\Pr \bigl( (D_i)_{i=0}^\ell = (d_i)_{i=0}^\ell \bigr)}
	            {\log_b(1 + b^{-\ell}/\bs{d}_{(\ell)})} - 1 \right|
	\ \le \ 3 \cdot h \bigl( \lfloor 36 \sigma^2 \rfloor \bigr) ,
$$
where $\bs{d}_{(\ell)} := \sum_{i=1}^\ell d_i \cdot b^{-i}$.	\hfill	$\Box$
\end{Corollary}

\subsection{Gumbel and Weibull distributions}

Let $X > 0$ be a random variable with Weibull distribution, i.e.\ for some parameters $\gamma, \tau > 0$,
$$
	\Pr(X \le r) \ = \ 1 - \exp( - (r/\gamma)^\tau )	\quad	\text{for} \ r \ge 0 .
$$
Then the standardized random variable $Y_o := \tau \log(X/\gamma)$ satisfies
$$
	F_o(y) \ := \ \Pr(Y_o \le y)
	\ = \ 1 - \exp( - e^y)	\quad\text{for} \ y \in \R
$$
and has density function
$$
	f_o(y) \ = \ e^y \exp(- e^y) ,
$$
i.e.\ $- Y_o$ has a Gumbel distribution. Thus $Y := \log_b(X)$ may be written as $Y = \mu + \sigma Y_o$ with $\mu := \log_b(\gamma)$ and $\sigma = (\tau \log b)^{-1}$.

Elementary calculations reveal that for any integer $n \ge 1$,
$$
	f_o^{(n-1)}(y) \ = \ p_n^{}(e^y) \exp(- e^y)
$$
with $p_n(t)$ being a polynomial in $t$ of degree $n$. Precisely, $p_1(t) = t$, and
\begin{equation}
	p_{n+1}^{}(t) \ = \ t (p_n'(t) - p_n^{}(t))
	\label{recursion}
\end{equation}
for $n = 1, 2, 3, \ldots$. In particular, $p_2(t) = t(1 - t)$ and $p_3(t) = t(1 - 3t + t^2)$. These considerations lead already to the following conclusion:

\begin{Corollary}
\label{cor: Gumbel-Weibull}
Let $X > 0$ have Weibull distribution with parameters $\gamma, \tau > 0$ as above. Then $\TV(f_o^{(k)}) < \infty$ and
$$
	\left| \frac{\Pr \bigl( (D_i)_{i=0}^\ell = (d_i)_{i=0}^\ell \bigr)}
	            {\log_b(1 + b^{-\ell}/\bs{d}_{(\ell)})} - 1 \right|
	\ \le \ 3 \cdot \TV(f_o^{(k)}) \Bigl( \frac{\tau \log b}{6} \Bigr)^{k+1}
$$
for arbitrary integers $k, \ell \ge 0$ and digits $d_0, d_1, d_2 \ldots$ as in Corollary~\ref{cor: Benford-Gauss}.	\hfill	$\Box$
\end{Corollary}

Explicit inequalities as in the gaussian case seem to be out of reach. Nevertheless some numerical bounds can be obtained. Table~\ref{tab: Gumbel-Weibull} contains numerical approximations for $\TV(f_o^{(k)})$ and the resulting upper bounds
$$
	B_{\tau}(k) \ := \ 3 \cdot \TV(f_o^{(k)}) \Bigl( \frac{\tau \, \log(10)}{6} \Bigr)^{k+1}
$$
for the maximal relative approximation error in Benford's law with decimal expansions, where $\tau = 1.0, 0.5, 0.3$. Note that $\tau = 1.0$ corresponds to the standard exponential distribution. For a detailed analysis of this special case we refer to Engel and Leuenberger (2003) and Miller and Nigrini (2006).

\begin{table}[h]
$$
	\begin{array}{|c||l|l|l|l|}
\hline
	 k & \mathrm{TV}(f_o^{(k)}) & B_{1.0}(k) & B_{0.5}(k) & B_{0.3}(k)
\\\hline\hline
	 0 & 7.3576 \cdot 10^{-1} & 8.4707 \cdot 10^{-1} & 4.2354 \cdot 10^{-1} & 2.5412 \cdot 10^{-1}
\\\hline
	 1 & 9.4025 \cdot 10^{-1} & 4.1543 \cdot 10^{-1} & 1.0386 \cdot 10^{-1} & 3.7388 \cdot 10^{-2}
\\\hline
	 2 & 1.7830               & 3.0232 \cdot 10^{-1} & 3.7790 \cdot 10^{-2} & 8.1627 \cdot 10^{-3}
\\\hline
	 3 & 4.5103               & \mathbf{2.9348 \cdot 10^{-1}}
	                                                 & 1.8343 \cdot 10^{-2} & 2.3772 \cdot 10^{-3}
\\\hline
	 4 & 1.4278 \cdot 10      & 3.5653 \cdot 10^{-1} & 1.1142 \cdot 10^{-2} & 8.6638 \cdot 10^{-4}
\\\hline
	 5 & 5.4301 \cdot 10      & 5.2038 \cdot 10^{-1} & 8.1309 \cdot 10^{-3} & 3.7936 \cdot 10^{-4}
\\\hline
	 6 & 2.4118 \cdot 10^2    & 8.8699 \cdot 10^{-1} & 6.9296 \cdot 10^{-3} & 1.9399 \cdot 10^{-4}
\\\hline
	 7 & 1.2252 \cdot 10^3    & 1.7292               & \mathbf{6.7546 \cdot 10^{-3}}
	                                                                        & 1.1345 \cdot 10^{-4}
\\\hline
	 8 & 7.0056 \cdot 10^3    & 3.7944               & 7.4110 \cdot 10^{-3} & 7.4686 \cdot 10^{-5}
\\\hline
	 9 & 4.4527 \cdot 10^4    & 9.2552               & 9.0383 \cdot 10^{-3} & 5.4651 \cdot 10^{-5}
\\\hline
	10 & 3.1140 \cdot 10^5    & 2.4840 \cdot 10      & 1.2129 \cdot 10^{-2} & 4.4003 \cdot 10^{-5}
\\\hline
	11 & 2.3763 \cdot 10^6    & 7.2744 \cdot 10      & 1.7760 \cdot 10^{-2} & 3.8659 \cdot 10^{-5}
\\\hline
	12 & 1.9648 \cdot 10^7    & 2.3083 \cdot 10^2    & 2.8177 \cdot 10^{-2} & \mathbf{3.6801 \cdot 10^{-5}}
\\\hline
	13 & 1.7498 \cdot 10^8    & 7.8888 \cdot 10^2    & 4.8150 \cdot 10^{-2} & 3.7732 \cdot 10^{-5}
\\\hline
	14 & 1.6698 \cdot 10^9    & 2.8890 \cdot 10^3    & 8.8166 \cdot 10^{-2} & 4.1454 \cdot 10^{-5}
\\\hline
	\end{array}
$$
\caption{Some bounds for Weibull-distributed $X$ with $\tau \le 1.0, 0.5, 0.3$}
\label{tab: Gumbel-Weibull}
\end{table}

\paragraph{Remark.}
Writing
$$
	p_n^{}(t) \ = \ \sum_{k=1}^n (-1)_{}^{k-1} S_{n,k}^{} \ t_{}^k ,
$$
it follows from the recursion (\ref{recursion}) that the coefficients can be calculated inductively by
\begin{equation*}
	S_{1,1}^{} \ = \ 1 ,	\qquad	S_{n,k}^{} = S_{n-1,k-1}^{} + k S_{n-1,k}^{} .
\end{equation*}
Hence the $S_{n,k}^{}$ are Stirling numbers of the second kind (see \cite{Graham_1994}, chapter 6.1).

\section{Proofs}

\subsection{Some useful facts about total variation}

In our proofs we shall utilize the some basic properties of total variation of functions $h : \mathbb{J} \to \R$ (cf.\ Royden 1988, Chapter~5). Note first that
$$
	\TV(h, \mathbb{J}) \ = \ \TV^+(h, \mathbb{J}) + \TV^-(h, \mathbb{J})
$$
with
$$
	\TV^{\pm}(h,\mathbb{J}) \ := \ \sup \Bigl\{ \sum_{i=1}^m \bigl( h(t_i) - h(t_{i-1}) \bigr)^{\pm} \, : \
		m \in \N; \, t_0 < t_1 < \cdots < t_m; \, t_0,\ldots,t_m \in \mathbb{J} \Bigr\}
$$
and $a^\pm := \max(\pm a,0)$ for real numbers $a$. Here are further useful facts in case of $\mathbb{J} = \R$:

\begin{Lemma}
\label{lem: TV and limits}
Let $h : \R \to \R$ with $\TV(h) < \infty$. Then both limits $h(\pm \infty) := \lim_{x \to \pm\infty} h(x)$ exist. Moreover, for arbitrary $x \in \R$,
$$
	h(x) \ = \ h(-\infty) + \TV^+(h, (-\infty,x]) - \TV^-(h, (-\infty,x]) .
$$
In particular, if $h(\pm \infty) = 0$, then $\TV^+(h) = \TV^-(h) = \TV(h)/2$.	\hfill	$\Box$
\end{Lemma}

\begin{Lemma}
\label{lem: TV and integral}
Let $h$ be integrable over $\R$.

\noindent
\textbf{(a)} \ If $\TV(h) < \infty$, then $\lim_{|x| \to \infty} h(x) = 0$.

\noindent
\textbf{(b)} \ If $h$ is $k\ge 1$ times absolutely continuous with $\TV(h^{(k)}) < \infty$ for some version of $h^{(k)}$, then
$$
	\lim_{|x| \to \infty} h^{(j)}(x) \ = \ 0
	\quad\text{for} \ j=0,1,\ldots,k .
$$
\end{Lemma}

While Lemma~\ref{lem: TV and limits} is standard, we provide a proof of Lemma~\ref{lem: TV and integral}:

\paragraph{Proof of Lemma~\ref{lem: TV and integral}.} Part~(a) follows directly from Lemma~\ref{lem: TV and limits}. Since $\TV(h) < \infty$, there exist both limits $\lim_{x \to \pm\infty} h(x)$. If one of these limits was nonzero, the function $h$ could not be integrable over $\R$.

For the proof of part~(b), define $h^{(k)}(\pm \infty) := \lim_{x \to \pm\infty} h^{(k)}(x)$. If $h^{(k)}(+\infty) \ne 0$, then one can show inductively for $j = k-1, k-2, \ldots, 0$ that $\lim_{x \to \infty} h^{(j)}(x) = \mathrm{sign}(h^{(k)}(+\infty)) \cdot \infty$. Similarly, if $h^{(k)}(-\infty) \ne 0$, then $\lim_{x \to -\infty} h^{(j)}(x) = (-1)^{k-j} \mathrm{sign}(h^{(k)}(-\infty)) \cdot \infty$ for $0 \le j < k$. In both cases we would get a contradiction to $h^{(0)} = h$ being integrable over $\R$.

Now suppose that $\lim_{|x| \to \infty} h^{(k)}(x) = 0$. It follows from Taylor's formula that for $x \in \R$ and $u \in [-1,1]$,
\begin{eqnarray*}
	|h(x+u)|
	& = & \left| \sum_{j=0}^{k-1} \frac{h^{(j)}(x)}{j!} \, u^j
		+ \int_0^u \frac{h^{(k)}(x + v) (u - v)^{k-1}}{(k-1)!} \, dv \right| \\
	& \ge & \Bigl| \sum_{j=0}^{k-1} \frac{h^{(j)}(x)}{j!} \, u^j \Bigr|
		- \sup_{|s| \ge |x|-1} \frac{|h^{(k)}(s)| |u|^k}{k!} .
\end{eqnarray*}
Hence
$$
	\int_{x-1}^{x+1} |h(t)| \, dt
	\ \ge \ \frac{|h^{(j)}(x)|}{j!} \, A_{j,k-1} - 2 \sup_{|s| \ge |x|-1} \frac{|h^{(k)}(s)|}{(k+1)!}
$$
for any $j \in \{0,1,\ldots,k-1\}$, where for $0 \le \ell \le m$,
$$
	A_{\ell,m} \ := \ \min_{a_0, \ldots, a_m \in \R \, : \, a_\ell = 1}
		\int_{-1}^1 \Bigl| \sum_{j=0}^m a_j u^j \Bigr| \, du
	\ > \ 0 .
$$
This shows that
$$
	|h^{(j)}(x)|
	\ \le \ \frac{j!}{A_{j,k-1}} \Bigl( \int_{x-1}^{x+1} |h(t)| \, dt
		+ 2 \sup_{|s| \ge |x| - 1} \frac{|h^{(k)}(s)|}{(k+1)!} \Bigr)
	\ \to \ 0	\quad\text{as} \ |x| \to \infty .
	\eqno{\Box}
$$

\subsection{Proofs of the main results}

\paragraph{Proof of Theorem~\ref{thm: TV(f0)}.}
For arbitrary $m \in \N$ and $0 \le t_0 < t_1 < \ldots < t_m \le 1$,
\begin{equation}
	\sum_{z \in \Z} \sum_{i=1}^m \bigl| f(z+t_i) - f(z+t_{i-1}) \bigr| \ \le \ \TV(f) .
	\label{eq: TV gf}
\end{equation}
In particular, for two points $x,y \in [0,1]$ with $\min(g(x), g(y)) < \infty$, the difference $g(x) - g(y)$ is finite. Hence $g < \infty$ everywhere. Now it follows directly from (\ref{eq: TV gf}) that $\TV(g) \le \TV(f)$. Moreover, for $0 \le x < y \le 1$,
\begin{eqnarray*}
	\bigl( g(y) - g(x) \bigr)^\pm
	& = & \Bigl( \sum_{z \in \Z} \bigl( f(z+y) - f(z+x) \bigr) \Bigr)^\pm \\
	& \le & \sum_{z \in \Z} \bigl( f(z+y) - f(z+x) \bigr)^\pm \\
	& \le & \TV^\pm(f) \\
	& = & \TV(f)/2 ,
\end{eqnarray*}
where the latter equality follows from Lemma~\ref{lem: TV and integral}~(a) and Lemma~\ref{lem: TV and limits}. 	\hfill	$\Box$

\paragraph{Proof of Corollary~\ref{cor: TV(f0)}.}
Let $0 \le x < y \le 1$ and $\delta := y - x \in (0,1]$. Then
\begin{eqnarray*}
	\bigl| G(y) - G(x) - (y - x) \bigr|
	& = & \Bigl| \int_x^y g(u) \, du - \delta \int_{y-1}^y g(u) \, du \Bigr| \\
	& = & \Bigl| (1 - \delta) \int_x^y g(u) \, du - \delta \int_{y-1}^{x} g(u) \, du \Bigr| \\
	& = & \Bigl| \delta (1 - \delta) \int_0^1 \bigl( g(x + \delta t) - g(x - (1 - \delta) t) \bigr) \, dt \Bigr| \\
	& \le & \delta (1 - \delta) \int_0^1 \bigl| g(x + \delta t) - g(x - (1 - \delta) t) \bigr| \, dt \\
	& \le & \delta (1 - \delta) \TV(f)/2 .	\qquad\qquad\qquad	\Box
\end{eqnarray*}

\paragraph{Proof of Theorem~\ref{thm: TV(fk)}.}
Throughout this proof let $x, y \in \R$ be generic real numbers with $\delta := y - x \in [0,1]$.
For integers $j \in \{0,\ldots,k\}$ and $N \ge 1$ we define
$$
	g_N^{(j)}(x,y) \ := \ \sum_{z=-N}^N \bigl( f^{(j)}(z+y) - f^{(j)}(z+x) \bigr) .
$$
Note that $g(y) - g(x) = \lim_{N \to \infty} g_N^{(0)}(x,y)$ whenever $g(x) < \infty$ or $g(y) < \infty$. To establish a relation between $g^{(j)}(\cdot,\cdot)$ and $g^{(j+1)}(\cdot,\cdot)$ note first that for absolutely continuous $h : \R \to \R$,
\begin{eqnarray*}
	h(y) - h(x)
	& = & h(y) - h(x) - \delta \bigl( h(y) - h(y-1) \bigr)
		+ \delta \bigl( h(y) - h(y-1) \bigr) \\
	& = & \delta (1 - \delta) \int_0^1 \bigl( h'(x + \delta t) - h'(x - (1 - \delta) t) \bigr) \, dt
		+ \delta \bigl( h(y) - h(y-1) \bigr) \\
	& = & \delta (1 - \delta) \int_0^1 \bigl( h'(x + \delta t) - h'(x + \delta t - t) \bigr) \, dt
		+ \delta \bigl( h(y) - h(y-1) \bigr) ,
\end{eqnarray*}
see also the proof of Corollary~\ref{cor: TV(f0)}. Hence for $0 < j \le k$,
\begin{eqnarray}
	g_N^{(j-1)}(x,y) & = & \delta (1 - \delta) \int_0^1 g_N^{(j)}(x + \delta t, x + \delta t - t) \, dt
	\label{eq: Recursion g_N} \\
	&& + \ \delta \bigl( f^{(j-1)}(N+y) - f^{(j-1)}(-N+y-1) \bigr) .
	\nonumber
\end{eqnarray}

Recall that $\lim_{|z| \to \infty} f^{(j)}(z) = 0$ for $0 \le j \le k$ by virtue of Lemma~\ref{lem: TV and integral}~(b). In particular, $\TV^\pm(f^{(k)}) = \TV(f^{(k)})/2$ by Lemma~\ref{lem: TV and limits}. Hence
$$
	g_N^{(k)}(x,y) \ = \ \sum_{z=-N}^N \bigl( f^{(k)}(y) - f^{(k)}(x) \bigr)^+
		- \sum_{z=-N}^N \bigl( f^{(k)}(y) - f^{(k)}(x) \bigr)^-
$$
satisfies the inequality $\left| g_N^{(k)}(x,y) \right| \le \TV(f^{(k)})/2$ and converges to a limit $g^{(k)}(x,y)$ as $N \to \infty$. Moreover, it follows from (\ref{eq: Recursion g_N}) that
$$
	\left| g_N^{(k-1)}(x,y) \right|
	\ \le \ \delta (1 - \delta) \TV(f^{(k)})/2 + 2 \|f^{(k-1)}\|_\infty
$$
and, via dominated convergence,
$$
	\lim_{N \to \infty} g_N^{(k-1)}(x,y)
	\ = \ g^{(k-1)}(x,y) := \delta (1 - \delta) \int_0^1 g^{(k)}(x + \delta t, x + \delta t - t) \, dt
$$
with
$$
	\left| g^{(k-1)}(x,y) \right|
	\ \le \ \delta (1 - \delta) \int_0^1 \left| g^{(k)}(x + \delta t, x + \delta t - t) \right| \, dt
	\ \le \ \delta (1 - \delta) \TV(f^{(k)})/2 .
$$

Now we perform an induction step: Suppose that for some $1 \le j < k$,
$$
	\left| g_N^{(j)}(x,y) \right| \ \le \ \alpha^{(j)} < \infty
$$
and
$$
	g^{(j)}(x,y) \ := \ \lim_{N \to \infty} g_N^{(j)}(x,y)
	\quad\text{exists with}\quad
	\left| g^{(j)}(x,y) \right| \ \le \ \beta^{(j)} \delta (1 - \delta) .
$$
For $j = k-1$ this is true with $\beta^{(k-1)} := \TV(f^{(k)})/2$. Now it follows from (\ref{eq: Recursion g_N}) and dominated convergence that
$$
	\left| g_N^{(j-1)}(x,y) \right| \ \le \ \alpha^{(j)} + 2 \|f^{(j-1)}\|_\infty
$$
and
$$
	\lim_{N \to \infty} g_N^{(j-1)}(x,y) \ = \ g^{(j-1)}(x,y)
	:= \delta (1 - \delta) \int_0^1 g^{(j)}(x + \delta t, x + \delta t - t) \, dt ,
$$
where
\begin{eqnarray*}
	\left| g^{(j-1)}(x,y) \right|
	& \le & \delta (1 - \delta) \int_0^1 \left| g^{(j)}(x + \delta t, x + \delta t - t) \right| \, dt \\
	& \le & \beta^{(j)} \delta (1 - \delta) \int_0^1 t(1 - t) \, dt \\
	& = & (\beta^{(j)}/6) \, \delta (1 - \delta) .
\end{eqnarray*}

These considerations show that $g^{(0)}(x,y) := \lim_{N \to \infty} g_N^{(0)}(x,y)$ always exists and satisfies the inequality
$$
	\left| g^{(0)}(x,y) \right|
	\ \le \ \delta (1 - \delta) \frac{\TV(f^{(k)})}{2 \cdot 6^{k-1}}
	\ \le \ \frac{\TV(f^{(k)})}{8 \cdot 6^{k-1}} .
$$
In particular, $g$ is everywhere finite with $g(y) - g(x) = g^{(0)}(x,y)$ satisfying the asserted inequalities.	\hfill	$\Box$

\paragraph{Proof of Corollary~\ref{cor: TV(fk)}.}
For $0 \le x < y \le 1$ and $\delta := y - x \in (0,1]$,
\begin{eqnarray*}
	\bigl| G(y) - G(x) - (y - x) \bigr|
	& = & \Bigl| \delta (1 - \delta)
		\int_0^1 \bigl( g(x + \delta t) - g(x + \delta t - t) \, dt \Bigr| \\
	& \le & \delta (1 - \delta) \frac{\TV(f^{(k)})}{2 \cdot 6^{k-1}} \int_0^1 t (1 - t) \, dt \\
	& = & \delta (1 - \delta) \frac{\TV(f^{(k)})}{2 \cdot 6^k} .
		\qquad\qquad\qquad	\Box
\end{eqnarray*}

\paragraph{Proof of Corollary~\ref{cor: f0}.} It is wellknown that integrable functions on the real line may be approximated arbitrarily well in $L^1(\R)$ by regular functions, for instance, functions with compact support and continuous derivative. With little extra effort one can show that for any fixed $\epsilon > 0$ there exists a probability density $\tilde{f}_o$ such that $\TV(\tilde{f}_o) < \infty$ and
$$
	\int_{-\infty}^\infty \bigl| f_o(z) - \tilde{f}_o(z) \bigr| \, dz \ \le \ \epsilon .
$$
With $\tilde{f}(x) := \tilde{f}_o((x - \mu)/\sigma)/\sigma$ and $\tilde{g}(x) := \sum_{z \in \Z} \tilde{f}(z+x)$,
$$
	\int_0^1 |g(x) - 1| \, dx
	\ \le \ \int_0^1 \bigl| g(x) - \tilde{g}(x) \bigr| \, dx + \int_0^1 |\tilde{g}(x) - 1| \, dx .
$$
But
\begin{eqnarray*}
	\int_0^1 \bigl| g(x) - \tilde{g}(x) \bigr| \, dx
	& \le & \int_0^1 \sum_{z \in \Z} \bigl| f(z+x) - \tilde{f}(z+x) \bigr| \, dx \\
	& = & \int_{-\infty}^\infty \bigl| f(y) - \tilde{f}(y) \bigr| \, dy \\
	& = & \int_{-\infty}^\infty \bigl| f_o(z) - \tilde{f}_o(z) \bigr| \, dz \\
	& \le & \epsilon
\end{eqnarray*}
while
$$
	\int_0^1 |\tilde{g}(x) - 1| \, dx
	\ \le \ \frac{\TV(\tilde{f})}{2}
	\ = \ \frac{\TV(\tilde{f}_o)}{2 \sigma}
	\ \to \ 0	\quad	(\sigma \to \infty)
$$
by means of Theorem~\ref{thm: TV(f0)}. Since $\epsilon > 0$ is arbitrarily small, this yields the asserted result.	\hfill	$\Box$

\paragraph{Proof of Theorem~\ref{thm: Benford-Gauss}.}
According to Theorem~\ref{thm: TV(f0)},
$$
	\mathrm{R}(g_{\mu,\sigma})
	\ \le \ \frac{\TV(f_{\mu,\sigma})}{2}
	\ = \ \frac{\TV(\phi)}{2\sigma}
	\ = \ \frac{\phi(0)}{\sigma} ,
$$
whereas Theorem~\ref{thm: TV(fk)} and the considerations in Section~\ref{subsec: Benford-Gauss} yield the inequalities
$$
	\mathrm{R}(g_{\mu,\sigma})
	\ \le \ \frac{\TV(f_{\mu,\sigma}^{(k)})}{8 \cdot 6^{k-1}}
	\ = \ \frac{\TV(\phi^{(k)})}{8 \cdot 6^{k-1} \sigma^{k+1}}
	\ \le \ \frac{\sqrt{ (k+1)! }}{8 \cdot 6^{k-1} \sigma^{k+1}}
$$
for all $k \ge 1$. Since the right hand side equals $0.75 / \sigma \ge \phi(0)/\sigma$ if we plug in $k = 0$, we may conclude that
$$
	\mathrm{R}(g_{\mu,\sigma})
	\ \le \ \frac{\sqrt{(k+1)!}}{8 \cdot 6^{k-1} \sigma^{k+1}}
	\ = \ 4.5 \cdot \sqrt{ \frac{(k+1)!}{(36 \sigma^2)^{k+1}} }
$$
for all $k \ge 0$. The latter bound becomes minimal if $k+1 = \lfloor 36 \sigma^2\rfloor \ge 1$, and this value yields the desired bound $4.5 \cdot h \bigl( \lfloor 36 \sigma^2\rfloor \bigr)$.

Similarly, Corollaries~\ref{cor: TV(f0)} and \ref{cor: TV(fk)} yield the inequalities
\begin{eqnarray*}
	\mathrm{KD}(G_{\mu,\sigma})
	& \le & \frac{\sqrt{(k+1)!}}{8 \cdot 6^k \sigma^{k+1}}
		\ = \ 0.75 \cdot \sqrt{ \frac{(k+1)!}{(36 \sigma^2)^{k+1}} } , \\
	\mathrm{MRAE}(G_{\mu,\sigma})
	& \le & \frac{\sqrt{(k+1)!}}{2 \cdot 6^k \sigma^{k+1}}
		\ = \ 3 \cdot \sqrt{ \frac{(k+1)!}{(36 \sigma^2)^{k+1}} } , \\
\end{eqnarray*}
for arbitrary $k \ge 0$, and $k+1 = \lfloor 36 \sigma^2\rfloor \ge 1$ leads to the desired bounds.	\hfill	 $\Box$

\paragraph{Acknowledgement.}
We are grateful to Steven J.\ Miller and an anonymous referee for constructive comments on previous versions of this manuscript.

\end{document}